%% file: efficient_FSBP_main.tex
\documentclass[final,onefignum,onetabnum,letterpaper]{siamonline250211}
\usepackage[scale=0.8]{geometry} 

\input{efficient_FSBP_shared}

\usepackage{lineno}
\usepackage{comment}

\ifpdf
\hypersetup{
	pdftitle={Gaussian FSBP operators},
 	pdfauthor={Glaubitz, Offner, and Stafforst}
}

\fi

\begin{document}

\maketitle

\begin{abstract}
	\input{0_abstract}
\end{abstract}





\input{1_intro} 
\input{2_FSBP}
\input{3_quadratures} 
\input{4_examples}
\input{5_numerics} 
\input{6_summary}

\section*{Acknowledgements}
\input{acknowledgements}


\bibliographystyle{siamplain}
\bibliography{references}

\end{document}

%% file: efficient_FSBP_shared.tex
\usepackage{soul}
\usepackage{lipsum}
\usepackage{amsfonts}
\usepackage{epstopdf}
\ifpdf
  \DeclareGraphicsExtensions{.eps,.pdf,.png,.jpg}
\else
  \DeclareGraphicsExtensions{.eps}
\fi

\usepackage{datetime}
\newdateformat{monthyeardate}{%
	\monthname[\THEMONTH] \THEDAY, \THEYEAR}

\usepackage{academicons}
\usepackage{xcolor}
\renewcommand{\orcid}[1]{\href{https://orcid.org/#1}{\textcolor[HTML]{A6CE39}{orcid.org/#1}}}

\usepackage{amsmath}
\allowdisplaybreaks
\usepackage{amssymb}
\usepackage{commath}
\usepackage{mathtools}
\usepackage{bbm}
\usepackage{bm}

\usepackage{color}
\usepackage{graphicx}
\usepackage[small]{caption}
\usepackage{subcaption}

\usepackage{relsize}
\usepackage{adjustbox}
\usepackage{algorithm}
\usepackage[noend]{algpseudocode}
\algnewcommand\algorithmicinput{\textbf{Input:}}
\algnewcommand\Input{\item[\algorithmicinput]}
\algnewcommand\algorithmicoutput{\textbf{Output:}}
\algnewcommand\Output{\item[\algorithmicoutput]}
\usepackage{booktabs}
\usepackage{siunitx}
\usepackage{tikz}
\usepackage{cancel}
\usepackage{pifont}
\usetikzlibrary{bayesnet} 

\usepackage{enumitem}
\setlist[enumerate]{leftmargin=.5in}
\setlist[itemize]{leftmargin=.5in}

\newsiamremark{remark}{Remark}
\newsiamremark{example}{Example}
\newsiamremark{hypothesis}{Hypothesis}
\crefname{hypothesis}{Hypothesis}{Hypotheses}
\newsiamthm{claim}{Claim}

\headers{Gaussian FSBP operators: Comparison and application}{J.\ Glaubitz, H.\ Haase, P.\ \"Offner, and F.\ Stafforst}

\title{Gaussian FSBP operators: Comparison and application to numerical methods for hyperbolic conservation laws\thanks{\monthyeardate\today 
\corresponding{Finnja Stafforst} 
}}
\author{ 
Jan Glaubitz\thanks{
Department of Mathematics, Link\"oping University, Sweden
(\email{jan.glaubitz@liu.se}, \orcid{0000-0002-3434-5563}; \email{henry.haase@liu.se}, \orcid{0009-0008-4065-0426})
} 
\and 
Henry Haase\footnotemark[2]
\and 
Philipp \"Offner\thanks{
Institute of Mathematics, Clausthal University of Technology, Germany
(\email{mail@philippoeffner.de}, \orcid{0000-0002-1367-1917}; \email{finnja.stafforst@tu-clausthal.de}, \orcid{0009-0007-1511-1300})
} 
\and 
Finnja Stafforst\footnotemark[3]
}

\usepackage{amsopn}

\usepackage{comment}

\definecolor{blue1}{RGB}{0,76,153}
\definecolor{green1}{RGB}{1,109,57}
\definecolor{red1}{RGB}{153,0,0}
\definecolor{orange1}{RGB}{255,127,80}

\usepackage[normalem]{ulem}

\usepackage[normalem]{ulem}

\DeclarePairedDelimiter\ceil{\lceil}{\rceil}

\DeclareMathOperator{\diag}{diag}

\DeclareMathOperator*{\argmin}{arg\,min}

\newcommand{\scp}[2]{\left\langle{#1,\, #2}\right\rangle} 
\renewcommand{\d}{\mathrm{d}} 
\newcommand{\intd}{\, \mathrm{d}}

\newcommand{\R}{\mathbb{R}}



%% file: 0_abstract.tex
Function-space summation-by-parts (FSBP) operators enable conservative and energy-stable numerical methods for hyperbolic conservation laws based on general, non-polynomial approximation spaces. 
Recent works show that using generalized Gaussian quadrature significantly reduces the number of grid points required compared to existing constructions that have mostly focused on equidistant grids. 
In this paper, we compare open and closed FSBP operators constructed with generalized Gaussian quadratures and apply them to numerically solve hyperbolic conservation laws.  
Furthermore, to support open node distributions, we extend the FSBP framework by introducing function-space exact extrapolation operators and operationalize them in numerical schemes for solving hyperbolic conservation laws. 
Our numerical experiments include the one-dimensional linear advection, non-viscous Burgers, and compressible Euler equations of gas dynamics. 
We observe that applying FSBP operators in numerical schemes can improve efficiency and accuracy. 
Notably, we demonstrate these advantages in more challenging time-dependent settings compared to other recent works on Gaussian FSBP operators.

%% file: 1_intro.tex
\section{Introduction} 
\label{sec:intro} 

Summation-by-parts (SBP) operators are a powerful tool for numerically solving ordinary and partial differential equations (ODEs \& PDEs), particularly hyperbolic conservation laws.
SBP operators mimic the integration-by-parts property at the discrete level, which is essential for ensuring conservation and stability in numerical schemes. 
In fact, many discretizations for hyperbolic conservation laws can be formulated in terms of SBP operators. Examples include finite difference (FD) \cite{kreiss1974finite,strand1994summation,mattsson2004summation},
finite volume \cite{nordstrom2001finite,nordstrom2003finite}, continuous Galerkin \cite{hicken2016multidimensional,hicken2020entropy,abgrall2020analysisI},
discontinuous Galerkin (DG) \cite{gassner2013skew,carpenter2014entropy,chan2018discretely}, 
radial basis function (RBF) \cite{glaubitz2024energy},
flux reconstruction \cite{ranocha2016summation,ranocha2018stability,offner2018stability}, 
active flux/PamPa \cite{barsukow2025stability,abgrall2025some}, 
and overset-grid \cite{glaubitz2025towards} methods. See also the review articles \cite{delreyfernandez2014review,svard2014review} and the more recent monographs \cite{ranocha2018generalised,glaubitz2020shock,oeffner2023approximation}.

\subsection*{Non-polynomial function spaces and FSBP operators}

Traditionally, SBP operators have been designed to be exact for polynomials up to a specific degree, assuming that these provide an accurate approximation to the solution of the PDE under consideration. 
However, in some cases, other function spaces provide more accurate approximations. 
Previous studies that highlighted the benefits of non-polynomial approximation spaces include exponentially fitted schemes for singular perturbation problems \cite{kadalbajoo2003exponentially,kalashnikova2009discontinuous}, DG methods \cite{yuan2006discontinuous} and (W)ENO reconstructions \cite{christofi1996study,iske1996structure,hesthaven2019entropy} based on non-polynomial approximation spaces, RBF schemes \cite{fornberg2015solving,fornberg2015primer}, and methods based on neural network approximations \cite{franck2024approximately}.
A key advantage of using non-polynomial approximation spaces is that one can leverage prior knowledge of the unknown solution's behavior. 

A framework of function-space SBP (FSBP) operators, which are SBP operators designed for general function spaces---beyond just polynomials---has recently been developed in a series of papers \cite{glaubitz2023summation,glaubitz2023multi,glaubitz2024summation,glaubitz2024energy,glaubitz2025optimization,glaubitz2025towards,glaubitz2026summation}.
These developments include first- and second-derivative FSBP operators \cite{glaubitz2023summation,glaubitz2024summation}, multi-dimensional FSBP operators \cite{glaubitz2023multi}, FSBP operators for global RBF methods \cite{glaubitz2024energy}, and sub-cell FSBP operators for overset grid methods \cite{glaubitz2025towards}. 
Over the last few years, the FSBP framework has received increased attention, as it enables the systematic investigation of the conservation and stability properties of existing non-polynomial PDE solvers and the development of new solver classes.

\subsection*{FSBP operators and quadratures}

The existence and construction of FSBP operators are intimately connected to classical quadrature theory. 
Let $\mathcal{F} \subset C^1$ be a finite-dimensional function space.
It has been shown in \cite[Corollary 4.6]{glaubitz2023summation} that there exists a $\mathcal{F}$-exact diagonal-norm FSBP operator if and only if there exists a positive and $(\mathcal{F}^2)'$-exact quadrature, where $(\mathcal{F}^2)' = \{ (fg)' \mid f,g \in \mathcal{F} \}$.
We will provide more details on FSBP operators and the required quadratures in \Cref{sec:FSBP,sec:quadratures}.

The first generation of FSBP operators \cite{glaubitz2023summation,glaubitz2023multi,glaubitz2024summation,glaubitz2024energy} has been developed for equidistantly and randomly distributed points and using least-squares quadratures \cite{huybrechs2009stable,glaubitz2020stableQF,glaubitz2021stableCFs,glaubitz2023construction}. 
The motivation for using these least-squares quadratures was that they can always be positive and $(\mathcal{F}^2)'$-exact as long as enough quadrature points are used, which was proven in \cite{glaubitz2023construction} under mild conditions on the point distribution. 
The resulting least-squares quadratures, however, use a number of points $N$ that is relatively large compared to the dimension of $\mathcal{F}$.\footnote{
It was numerically observed in several works \cite{huybrechs2009stable,glaubitz2020stableDG,glaubitz2020stableQF,glaubitz2021stableCFs,glaubitz2023construction} that $N \sim \dim((\mathcal{F}^2)')^2$ for equidistant points. 
} 
Preferably, the number of quadrature points should be the same as the dimension of $\mathcal{F}$, i.e., $N = \dim(\mathcal{F})$, resulting in so-called \emph{interpolatory FSBP operators}. 
A smaller number of points is not possible, and a larger number of points results in (i) a less-efficient FSBP operator and (ii) unresolved modes, which can result in (i) increased computation costs and (ii) increased numerical errors.  
See \cite{glaubitz2026summation} for more details on the disadvantages of having unresolved modes in FSBP operators. 

More recently, \cite{glaubitz2025optimization,glaubitz2026summation} explored optimization approaches for constructing FSBP operators in tandem with a required positive and $(\mathcal{F}^2)'$-exact quadrature. 
While such optimization approaches can result in more efficient FSBP operators, as they might result in quadratures with fewer points than the least squares quadrature, the resulting FSBP operators from \cite{glaubitz2025optimization,glaubitz2026summation} were still constructed mainly on equidistant points and typically used a number of points $N$ that is relatively large compared to the dimension of $\mathcal{F}$.

\subsection*{Our contribution: Leveraging generalized Gaussian quadratures}

Our goal is to construct efficient FSBP operators whose number of grid points $N$ is as close as possible to the dimension $\dim(\mathcal{F})$ of the underlying function space $\mathcal{F}$ and \emph{operationalize} them in numerical schemes for solving hyperbolic conservation laws.
Ideally, one would like to achieve $N = \dim(\mathcal{F})$, resulting in an interpolatory FSBP operator.
For polynomial function spaces $\mathcal{F}=\mathcal{P}_d$, where $\mathcal{P}_d$ denotes the space of polynomials of degree at most $d$, such optimal constructions are provided by Gauss--Lobatto and Gauss--Legendre quadratures. 
The Gauss--Lobatto quadrature is a \emph{closed} since it includes both interval endpoints, and achieves exactness for $\mathcal{P}_{2N-3}$ using $N$ nodes.
In contrast, the Gauss--Legendre quadrature is an \emph{open} quadrature rule, since it does not include the interval endpoints, and achieves exactness for $\mathcal{P}_{2N-1}$ using only $N$ nodes. 
It was demonstrated in \cite{offner2019error} that SBP operators based on such open Gauss--Legendre quadratures can yield numerical PDE solutions with higher accuracy for the same degrees of freedom compared to those based on closed Gauss--Lobatto quadratures; also see \cite{gassner2011comparison}. 
Generalizations of these quadrature rules to non-polynomial function spaces are provided by \emph{open} and \emph{closed generalized Gaussian quadratures} (GGQs); see \cite{ma1996generalized,cheng1999nonlinear,huybrechs2009generalized,bremer2010nonlinear,huybrechs2022computation}. 

In this work, we investigate the use of GGQs to construct efficient FSBP operators and operationalize them in numerical methods to solve hyperbolic conservation laws.
In particular, we study both closed FSBP operators based on closed GGQs and open FSBP operators based on open GGQs, and we compare their respective efficiencies. 
The recent work \cite{hale2026summation} established that closed GGQs yield more efficient closed FSBP operators than using (nearly) equidistant points together with a least-squares ansatz, extending the classical connection between Gauss--Lobatto quadratures and polynomial SBP operators. 
The use of open GGQs to construct FSBP operators has \textbf{not} yet been explored in depth.
The main contribution of this paper is to close this gap.\footnote{
During the final stages of this work, a related study \cite{bercik2026construction} introduced a construction of open and half-open FSBP operators based on the framework of Huybrechs \cite{huybrechs2022computation}. 
We were unaware of this recent work until its publication on arXiv on July 13; 
Our work had begun earlier, as demonstrated by various presentations at different conferences, including a workshop in Poland in March 2026 and the ICCFD13 (F. Stafforst: Function-Space SBP Operators for Structure-Preserving Discretizations, 8th of July). 
A detailed comparison of the adapted construction process and its application to solving hyperbolic problems is now planned for future work. 
} 
We show that open GGQs provide a natural foundation for open FSBP operators and, by avoiding the need to include boundary nodes, yield operators with fewer grid points than their closed-GGQ counterparts. 
Consequently, while closed GGQs remain optimal within the class of closed FSBP operators, open GGQs yield even more efficient FSBP discretizations whenever an open formulation is admissible. 
In this sense, our results extend the notion of optimality beyond the closed setting considered in \cite{hale2026summation} and align with the recent construction framework presented in \cite{bercik2026construction}.

Moreover, we operationalize and compare the resulting closed and open Gaussian FSBP operators for several hyperbolic conservation laws, which we believe to be significantly more challenging compared to the numerical experiments for FSBP operators presented in \cite{hale2026summation,bercik2026construction}.
Our numerical experiments include the one-dimensional linear advection, non-viscous Burgers, and compressible Euler equations of gas dynamics. 
In contrast, \cite{hale2026summation} considered a linear advection and a linear advection-diffusion problem, while \cite{bercik2026construction} considered three time-independent problems of the form $\frac{\d u}{\d x} = f(x,u)$.

At the same time, the accuracy gains of open FSBP operators come at a price. 
Because no grid points coincide with the interval endpoints, boundary and inter-element coupling conditions cannot be imposed directly at grid nodes; 
Instead, the numerical solution and its fluxes must be extrapolated to the boundary \cite{fernandez2014generalized}. 
In the FSBP setting, these extrapolation operators must be $\mathcal{F}$-exact, adding a construction step absent for closed FSBP operators. 
To this end, we extend the definition and construction of FSBP operators to grid points that exclude the endpoints of the interval. 
This extension has been developed for polynomial SBP operators in \cite{fernandez2014generalized}; see also \cite{ranocha20218general,chan2019efficient,offner2019error}. 
More recently, a related construction in the same spirit was presented in \cite{bercik2026construction}.

\subsection*{Outline}

The remainder of this paper is organized as follows. 
In \Cref{sec:FSBP}, we extend the FSBP framework to generalized FSBP operators whose grids do not necessarily include the boundary points, including the construction of $\mathcal{F}$-exact extrapolation and boundary operators. 
\Cref{sec:quadratures} addresses the required positive and $(\mathcal{F}^2)'$-exact quadratures: we recall open and closed GGQs, discuss their existence for Chebyshev sets, and describe the algorithm used to compute them. 
In \Cref{sec:examples}, we present examples of open Gaussian FSBP operators, and compare them with existing least-squares-based operators on equidistant points and with closed Gaussian FSBP operators. 
\Cref{sec:numerics} operationalizes the proposed operators in multi-block FSBP-SAT and FSBP-DG schemes and reports numerical experiments for hyperbolic conservation laws, comparing accuracy and efficiency of open, closed, and equidistant-grid FSBP operators. 
Finally, \Cref{sec:summary} offers concluding remarks and directions for future work.

%% file: 2_FSBP.tex
\section{Generalized FSBP operators} 
\label{sec:FSBP} 

We present a framework for generalized FSBP operators that do not necessarily include boundary points in their grids. 
This will subsequently allow us to construct highly efficient—and, in many cases, even optimal—versions of FSBP operators based on GGQs.

\subsection{Generalized FSBP operators}
\label{sub:FSBP}

Let $\Omega = [x_L,x_R] \subset \R$ be a compact interval and $\mathbf{x} = [x_1,\dots,x_N]$ be a vector of grid points on $\Omega$, where $x_L \leq x_1 < \dots < x_N \leq x_R$.
The existing literature on FSBP operators assumes that the grid includes the boundary points of the computational domain, i.e., $x_1 = x_L$ and $x_N = x_R$. 
However, this assumption fails for the quadrature points of the open GGQs we use here to construct more efficient FSBP operators.\footnote{The same is true for polynomial SBP operators when, for instance, the (open) Gauss--Legendre quadrature is used instead of the (closed) Gauss--Lobatto quadrature \cite{fernandez2014generalized}.}
We therefore focus on generalized FSBP operators that do not necessarily include one or both boundary nodes. Generalized SBP operators for polynomial function spaces were first introduced in \cite{fernandez2014generalized} (also see \cite{ranocha20218general,chan2019efficient,offner2019error}). \\
Furthermore, \cite{glaubitz2025towards} considered generalized SBP operators in the context of sub-cell SBP operators for overset grid methods. 
Here, we extend the concepts of generalized SBP operators to FSBP operators for non-polynomial function spaces.

\begin{definition}[Generalized FSBP operators]
\label{def:FSBP} 
	Let $\mathcal{F} \subset C^1([x_L,x_R])$ be a finite-dimensional function space. 
	An operator $D = P^{-1} Q \in \R^{N \times N}$ approximating $\partial_x$ is called an \emph{$\mathcal{F}$-exact generalized FSBP operator} if 
	\begin{enumerate} 
		\item[(i)] 
		$D \mathbf{f} =\mathbf{f}' $ for all $f \in \mathcal{F}$, 	
			
		\item[(ii)] 
		$P = \diag(\mathbf{p})$ with $\mathbf{p} = [p_1,\dots,p_N]$ and $p_n > 0$ for all $n=1,\dots,N$,  
		
		\item[(iii)] 
		 $Q + Q^T = B$, and  
		 
		\item[(iv)] 
		the boundary matrix $B \in \R^{N \times N}$ satisfies 
		\begin{equation}\label{eq:boundary_matrix_1d}
			\mathbf{f}^T B \mathbf{g} = fg|_{x_L}^{x_R}, \quad \forall f,g \in \mathcal{F}. 
		\end{equation}
			
	\end{enumerate}  
\end{definition} 

Relation (i) ensures that $D$ accurately approximates $\partial_x$ by requiring $D$ to be exact for all functions from $\mathcal{F}$. 
Condition (ii) guarantees that $P$ induces a discrete inner product and norm, which are given by $\scp{\mathbf{u}}{\mathbf{v}}_P = \mathbf{u}^T P \mathbf{v}$ and $\|\mathbf{u}\|^2_P = \mathbf{u}^T P \mathbf{u}$ for $\mathbf{u},\mathbf{v} \in \R^N$, respectively. 
Relation (iii) encodes the SBP property, which allows us to mimic integration-by-parts (IBP) on a discrete level. 
Recall that IBP reads 
\begin{equation}\label{eq:IBP_1d}
	\int_{x_L}^{x_R} u (\partial_x v) \intd x + \int_{x_L}^{x_R} (\partial_x u) v \intd x
		= u v \big|_{x_L}^{x_R} \quad \forall u, v \in C^1([x_L,x_R]).
\end{equation}
The discrete version of \cref{eq:IBP_1d}, which follows from (iii) in \cref{def:FSBP}, is 
\begin{equation}\label{eq:SBP_1d}
	\mathbf{u}^T P ( D \mathbf{v} ) + ( D \mathbf{u} )^T P \mathbf{v} = \mathbf{u}^T B \mathbf{v} \quad \forall \mathbf{u}, \mathbf{v} \in \R^N. 
\end{equation} 
Note that the two terms on the left-hand side of \cref{eq:SBP_1d} approximate the related terms on the left-hand side of \cref{eq:IBP_1d}.
Finally, (iv) in \cref{def:FSBP} ensures that also the right-hand side of \cref{eq:SBP_1d} accurately approximates the right-hand side of \cref{eq:IBP_1d}.  

A few remarks are in order. 

\begin{remark} 
	As mentioned above, existing FSBP operators have predominantly been constructed on equidistant (see \cite{glaubitz2023summation,glaubitz2024summation,glaubitz2024energy,glaubitz2026summation}) and random (see \cite{glaubitz2023multi,glaubitz2025optimization}) grids that include the boundary points of the computational domain. 
	In this case, $x_1 = x_L$, $x_N = x_R$ and $f_{1} = f(x_L)$, $f_{N} = f(x_R)$ for any function $f \in C^1([x_L,x_R])$. 
	Hence, $B = \diag(-1,0,\dots,0,1)$ satisfies $\mathbf{f}^T B \mathbf{g} = fg|_{x_L}^{x_R}$ for all $f,g \in C([x_L,x_R])$. 
	This shows that \cref{def:FSBP} includes these existing FSBP operators as a special case.  
\end{remark}

\begin{remark} 
	The extension of polynomial SBP operators to grids that do not necessarily include one or both boundary nodes has been described, for instance, in \cite{fernandez2014generalized,ranocha20218general,chan2019efficient,mateo2023flux} and references therein.  
	The resulting operators are usually referred to as ``generalized SBP operators". 
	Here, we further expand this concept to cover general function spaces, moving beyond the confines of traditional polynomial operators. 
	In the context of the present work, this is necessary to make generalized FSBP applicable to RBF approximation spaces.  
\end{remark}

\subsection{Constructing generalized FSBP operators}
\label{sub:const_FSBP}

We next briefly review the existence and construction of generalized FSBP operators.  
The existence of $\mathcal{F}$-exact FSBP operators has been characterized in \cite{glaubitz2023summation} in terms of quadratures. 
Although \cite{glaubitz2023summation} assumed that the grid included the boundary points of $[x_L,x_R]$, it is straightforward to see that the arguments carry over to the generalized \cref{def:FSBP}. 
We summarize this in \cref{lem:diag_char}. 

\begin{lemma}\label{lem:diag_char} 	
	Let $\mathcal{F} \subset C^1([x_L,x_R])$ be a function space with basis $\{ f_k \}_{k=1}^K$. 
	Assume that the Vandermonde matrix $V = [ \mathbf{f_1}, \dots, \mathbf{f_K} ]$ has linearly independent columns and that there exists a boundary operator $B$ satisfying \cref{eq:boundary_matrix_1d}. 
	Then, there exists an $\mathcal{F}$-exact FSBP operator $D=P^{-1}Q$ with a positive definite diagonal-norm matrix $P=\diag(\mathbf{p})$ if and only if $\mathbf{x}$ and $\mathbf{p}$ are the points and weights of a positive and $(\mathcal{F}^2)'$-exact quadrature on $[x_L,x_R]$.  
\end{lemma}

By $\mathbf{x}$ and $\mathbf{p}$ being the points and weights of a positive and $(\mathcal{F}^2)'$-exact quadrature on $[x_L,x_R]$ we mean that the quadrature  $I_N[f] = \sum_{n=1}^N p_n f(x_n)$---with points $\mathbf{x}$ and weights $\mathbf{p}$---satisfies $p_n > 0$ for all $n=1,\dots,N$ and $I_N[f] = \int_{x_L}^{x_R} f(x) \intd x$ for all $f \in (\mathcal{F}^2)'$. 
Here, 
\begin{equation}\label{eq:F_sq_prime}
	\left( \mathcal{F}^2 \right)' = \{ \, f' g + f g' \mid f,g \in \mathcal{F} \, \} 
\end{equation}	
is the space of all functions corresponding to the derivative of the product of two functions from $\mathcal{F}$. 
For completeness, we briefly revisit the construction of $\mathcal{F}$-exact FSBP operators using the procedure established in \cite{glaubitz2023summation} (also see \cite{glaubitz2023multi} for the multi-dimensional case and \cite{fernandez2014generalized,hicken2016multidimensional} for polynomial SBP operators). 

First, we construct a positive quadrature that is exact for $( \mathcal{F}^2 )'$. 
We then choose the positive weights $\mathbf{p}$ of this quadrature as the diagonal elements of the norm matrix $P$, i.e., $P = \diag(\mathbf{p})$.
Next, we get $Q$ by re-writing it as $Q = Q_A + B/2$, where $Q_A$ is the anti-symmetric part of $Q$, and recovering $Q_A$ by solving 
\begin{equation}\label{eq:exactness_QA} 
	Q_A V = P V' - \frac{1}{2} B V.
\end{equation}
Here, $V = [\mathbf{f_1},\dots,\mathbf{f_K}]$ and $V' = [\mathbf{f_1'},\dots,\mathbf{f_K'}]$ are the Vandermonde matrices for the basis elements $f_1,\dots,f_L$ of $\mathcal{F}$ and their derivatives, respectively. 
One can solve the matrix equation \cref{eq:exactness_QA} by recasting it as a linear system 
\begin{equation}\label{eq:LS_q} 
 	C \mathbf{q} = \mathbf{y}, 
\end{equation}
where $\mathbf{q}$ denotes the vector that contains the strictly lower part of $Q_A$. 
It is computationally convenient to select the unique least-squares solution of \cref{eq:LS_q}, which is the solution with minimal Euclidean norm $\|\cdot\|_{\ell^2}$ among all possible solutions \cite{golub2012matrix}.
The strictly upper part of $Q_A$ is then obtained by $(Q_A)_{j,i} = - (Q_A)_{i,j}$, and the diagonal elements are set to zero. 
Finally, we get an $\mathcal{F}$-exact FSBP operator as $D = P^{-1} Q$. 
Observe that the above procedure does not assume that the grid points include boundary points;
It works for general boundary operators $B$.

\begin{remark}
	As an alternative to the construction procedure above, recent works \cite{glaubitz2025optimization,glaubitz2026summation} considered optimization-based approaches for constructing FSBP operators. 
	While this methodology is new for FSBP operators, similar optimization procedures have previously been employed in meshfree methods to construct polynomial SBP operators on point clouds \cite{jameson2012,hicken2025constructing}. 
\end{remark}

\subsection{Constructing the boundary operators}
\label{sub:const_B_FSBP}

It remains to show that a boundary operator $B$ satisfying \cref{eq:boundary_matrix_1d} exists. 
Then, we can use the above procedure to construct an $\mathcal{F}$-exact FSBP operator. 
For polynomial SBP operators, \cite{fernandez2014generalized} suggested constructing the boundary matrix $B$ based on the idea of using polynomial interpolants of degree $N-1$ on the grid $\mathbf{x} = [x_1,\dots,x_N]$. 
Here, we modify this idea in two ways: 
\begin{enumerate}
	\item 
	We use a least squares approximation instead of interpolation, as this is more stable;  

	\item 
	We consider function approximations from $\mathcal{F}$ instead of polynomial ones. 

\end{enumerate}

To elaborate on the least squares approach, we follow \cite[Chapter 3.1]{glaubitz2020shock} and consider the discrete inner product $\langle f, g \rangle_{\mathbf{x}} = \sum_{n=1}^N w_n f(x_n) g(x_n)$ for $f,g \in \mathcal{F}$. 
Then, for each $f \in C^{1}([x_L,x_R])$, there exists a unique least squares approximation $f_{\mathbf{x}} \in \mathcal{F}$ such that 
\begin{equation} 
	f_{\mathbf{x}} = \argmin_{g \in \mathcal{F}} \left\{ \| f - g \|_{\mathbf{x}} \right\},
\end{equation}
where $\| f \|_{\mathbf{x}}^2 = \langle f, f \rangle_{\mathbf{x}}$ is the induced norm. 
Evaluating the least squares approximation at the boundary points yields 
\begin{equation}\label{eq:evaluate_LS} 
	\mathbf{t}_L^T \mathbf{f} = f_{\mathbf{x}}(x_L) \approx f(x_L), \quad 
	\mathbf{t}_R^T \mathbf{f} = f_{\mathbf{x}}(x_R) \approx f(x_R),
\end{equation}
where it remains to address how the vectors $\mathbf{t}_L, \mathbf{t}_R \in \R^N$ look like. 
Assuming that \cref{eq:evaluate_LS} holds, we form the \emph{extrapolation matrix} 

\begin{equation} 
	E = \begin{bmatrix} 1 \\ 0 \end{bmatrix} \mathbf{t}_L^T + \begin{bmatrix} 0 \\ 1 \end{bmatrix} \mathbf{t}_R^T
\end{equation} 
and get a desired boundary operator as 
\begin{equation}\label{eq:LS_boundary_op}
	B = E^T \begin{bmatrix} -1 & 0 \\ 0 & 1 \end{bmatrix} E.
\end{equation} 

To see that $B$ in \cref{eq:LS_boundary_op} satisfies \cref{eq:boundary_matrix_1d}, recall that the least squares approximation is exact for all functions from $\mathcal{F}$, i.e., $f_{\mathbf{x}} = f$ for all $f \in \mathcal{F}$. 
Hence, $\mathbf{t}_L^T \mathbf{f} = f(x_L)$ and $\mathbf{t}_R^T \mathbf{f} = f(x_R)$ for all $f \in \mathcal{F}$. 
This, in turn, implies 

\begin{equation} 
\begin{aligned} 
	\mathbf{f}^T B \mathbf{g} 
		& = \left( E \mathbf{f} \right)^T \begin{bmatrix} -1 & 0 \\ 0 & 1 \end{bmatrix} \left( E \mathbf{g} \right) \\ 
		& = \left( \begin{bmatrix} 1 \\ 0 \end{bmatrix} f(x_L) + \begin{bmatrix} 0 \\ 1 \end{bmatrix} f(x_R) \right)^T \begin{bmatrix} -1 & 0 \\ 0 & 1 \end{bmatrix} \left( \begin{bmatrix} 1 \\ 0 \end{bmatrix} g(x_L) + \begin{bmatrix} 0 \\ 1 \end{bmatrix} g(x_R) \right) \\ 
		& = \begin{bmatrix} f(x_L) \\ f(x_R) \end{bmatrix}^T \begin{bmatrix} -1 & 0 \\ 0 & 1 \end{bmatrix} \begin{bmatrix} g(x_L) \\ g(x_R) \end{bmatrix} \\
		& = f(x_R) g(x_R)-f(x_L) g(x_L)
\end{aligned}
\end{equation} 
for all $f,g \in \mathcal{F}$,
which shows that \cref{eq:boundary_matrix_1d} holds. 

It now remains to derive explicit formulas for $\mathbf{t}_L$ and $\mathbf{t}_R$.
To this end, let $\{ f_k \}_{k=1}^K$ be an arbitrary basis of $\mathcal{F}$.
Then, the least squares approximation $f_{\mathbf{x}}$ can be expressed as 
\begin{equation} 
	f_{\mathbf{x}} = \gamma_1 f_1 + \dots + \gamma_K f_K
\end{equation}
with the coefficients $\boldsymbol{\gamma} = [\gamma_1,\dots,\gamma_K]^T$ satisfying 
\begin{equation} 
	G \boldsymbol{\gamma} = \mathbf{\hat{f}}.
\end{equation}
Here, the Gram matrix $G$ and the vector of modal coefficients $\mathbf{\hat{f}}$ are  
\begin{equation}
	G = 
	\begin{bmatrix} 
		\langle f_1, f_1 \rangle_{\mathbf{x}} & \dots & \langle f_1, f_K \rangle_{\mathbf{x}} \\ 
		\vdots & & \vdots \\ 
		\langle f_K, f_1 \rangle_{\mathbf{x}} & \dots & \langle f_K, f_K \rangle_{\mathbf{x}}
	\end{bmatrix}, \quad 
	\mathbf{\hat{f}} = 
	\begin{bmatrix} 
		\langle f_1, f \rangle_{\mathbf{x}} / \| f_1 \|_{\mathbf{x}}^2 \\ 
		\vdots \\ 
		\langle f_K, f \rangle_{\mathbf{x}} / \| f_K \|_{\mathbf{x}}^2 
	\end{bmatrix}.
\end{equation} 
Observe that $G = V^T V$. 
If we want to evaluate the least squares approximation $f_{\mathbf{x}}$ at the boundary points, we get these function values as 
\begin{equation}\label{eq:evaluate_LS2} 
\begin{aligned}
	f_{\mathbf{x}}(x_L) 
		& = \mathbf{v}_L^T \boldsymbol{\gamma} 
		= \mathbf{v}_L^T G^{-1} \mathbf{\hat{f}}, \quad 
		&& \mathbf{v}_L^T = [ f_1(x_L), \dots, f_K(x_L) ], \\ 
	f_{\mathbf{x}}(x_R) 
		& = \mathbf{v}_R^T \boldsymbol{\gamma} 
		= \mathbf{v}_R^T G^{-1} \mathbf{\hat{f}}, \quad 
		&& \mathbf{v}_R^T = [ f_1(x_R), \dots, f_K(x_R) ].
\end{aligned}
\end{equation}
To simplify the above expressions and to connect the modal and nodal coefficients, $\mathbf{\hat{f}}$ and $\mathbf{f}$, we now restrict the basis $\{ f_k \}_{k=1}^K$ to be orthogonal w.r.t.\ the discrete inner product $\langle \cdot, \cdot \rangle_{\mathbf{x}}$.
In this case, $\langle f_k, f_l \rangle_{\mathbf{x}} = \delta_{k,l}$ and $\hat{f}_k = \sum_{n=1}^N f_k(x_n) f(x_n)$. 
Hence, $G = I$ and $\mathbf{\hat{f}} = V^T \mathbf{f}$, and \cref{eq:evaluate_LS2} becomes 
\begin{equation}\label{eq:evaluate_LS3} 
	f_{\mathbf{x}}(x_L) = \mathbf{v}_L^T V^T \mathbf{f}, \quad 
	f_{\mathbf{x}}(x_R) = \mathbf{v}_R^T V^T \mathbf{f}.
\end{equation}
This shows that $\mathbf{t}_L^T \mathbf{f} = f_{\mathbf{x}}(x_L)$ and $\mathbf{t}_R^T \mathbf{f} = f_{\mathbf{x}}(x_R)$ are satisfied by choosing $\mathbf{t}_L^T = \mathbf{v}_L^T V^T$ and $\mathbf{t}_R^T = \mathbf{v}_R^T V^T$, respectively. 

We summarize the resulting construction of generalized $\mathcal{F}$-exact FSBP operators in \cref{alg:FSBP_construction}.

\begin{algorithm}[t]
\caption{Construction of a generalized $\mathcal{F}$-exact FSBP operator $D = P^{-1} Q$}
\label{alg:FSBP_construction}
\begin{algorithmic}[1]
	\Input Interval $[x_L,x_R]$ and function space $\mathcal{F} \subset C^1([x_L,x_R])$ with basis $\{ f_k \}_{k=1}^K$.
	\Output Diagonal-norm $\mathcal{F}$-exact FSBP operator $D = P^{-1} Q$ in the sense of \cref{def:FSBP}.
	\State Compute a positive and $(\mathcal{F}^2)'$-exact quadrature on $[x_L,x_R]$ with nodes $\mathbf{x} = [x_1,\dots,x_N]$ and weights $\mathbf{p} = [p_1,\dots,p_N]$ \Comment{open or closed GGQ; see \Cref{sec:quadratures}}
	\State Set $P = \diag(\mathbf{p})$
	\State Orthonormalize the basis $\{ f_k \}_{k=1}^K$ w.r.t.\ $\scp{f}{g}_{\mathbf{x}} = \sum_{n=1}^N f(x_n) g(x_n)$ (e.g., via a QR factorization), and assemble $V = [\mathbf{f_1},\dots,\mathbf{f_K}]$ and $V' = [\mathbf{f_1'},\dots,\mathbf{f_K'}]$, so that $V^T V = I$
	\State Evaluate the basis at the boundary, giving $\mathbf{v}_L = [f_1(x_L),\dots,f_K(x_L)]^T$ and $\mathbf{v}_R = [f_1(x_R),\dots,f_K(x_R)]^T$, and form the extrapolation vectors $\mathbf{t}_L = V \mathbf{v}_L$ and $\mathbf{t}_R = V \mathbf{v}_R$
	\State Assemble the boundary operator $B = \mathbf{t}_R \mathbf{t}_R^T - \mathbf{t}_L \mathbf{t}_L^T$
	\State Recover the antisymmetric part $Q_A$ from $Q_A V = P V' - \tfrac{1}{2} B V$ by recasting it as the linear system $C \mathbf{q} = \mathbf{y}$ for the strictly lower triangular entries $\mathbf{q}$ of $Q_A$ and taking the minimal-norm least-squares solution; set $(Q_A)_{j,i} = -(Q_A)_{i,j}$ and $(Q_A)_{i,i} = 0$
\State Set $Q = Q_A + \tfrac{1}{2} B$
\State \Return $D = P^{-1} Q$
\end{algorithmic}
\end{algorithm}

%% file: 3_quadratures.tex
\section{GGQs and their application to FSBP operators} 
\label{sec:quadratures} 

We can always use the least-squares approach \cite{wilson1970discrete,huybrechs2009stable, glaubitz2021stableCFs,glaubitz2023construction} to find quadratures for $(\mathcal{F}^2)'$. 
In certain situations, more efficient GGQs  \cite{ma1996generalized,bremer2010nonlinear,yarvin1998generalized,huybrechs2022computation} can also be used, however. 
Here, we focus on these GGQs because they require fewer grid points than least-squares quadratures on equidistant points to achieve the same level of accuracy in the same function space.

\begin{definition}[Open and closed GGQs]
	Consider an interval $[a,b]$ and a system of $2m$ linearly independent functions $\{ \phi_1, \dots, \phi_{2m} \}$. 
	An \emph{open GGQ} for this system is a set of $m$ nodes $a < x_1 < \dots < x_m < b$ and weights $w_1,\dots,w_m$ such that 
	\begin{equation}
    		\int_{a}^{b} \phi_i \intd x = \sum_{i=1}^m w_i \phi_i(x_i), \quad i=1,\dots,2m.
	\end{equation}
	In contrast, a \emph{closed GGQ} for the same system uses a set of $m+1$ nodes $a = x_1 < \dots < x_{m+1} = b$ (including the boundary points) and weights $w_1,\dots,w_{m+1}$ such that 
	\begin{equation}
    		\int_{a}^{b} \phi_i \intd x = \sum_{i=1}^m w_i \phi_i(x_i), \quad i=1,\dots,2m.
	\end{equation}
\end{definition}

Notably, open GGQs exclude the interval endpoints $a$ and $b$. 
Consequently, the values at these two boundary points are not prescribed \emph{a priori}, leaving two additional degrees of freedom. 
As a result, open GGQs require fewer quadrature points than closed GGQs to achieve exactness for the same function space. 
In the present context, the function space of interest is $(\mathcal{F}^2)'$, as a positive $(\mathcal{F}^2)'$-exact quadrature induces a $\mathcal{F}$-exact FSBP operator.

We therefore expect open GGQs to yield more efficient FSBP operators than closed GGQs, in the sense that fewer grid points are needed to attain exactness over the same function space.

\subsection{Existence of open GGQs for Chebyshev sets}
\label{sub:chebyshev_quad}

In \cite{huybrechs2022computation}, Huybrechs presents an algorithm for computing open GGQs on Chebyshev sets. 
We use and adapt his approach to construct open Gaussian FSBP operators due to its theoretical foundations. 

For completeness, we first introduce the necessary definitions and then outline the underlying algorithm. 
We first define a Chebyshev set and describe the moment space.\footnote{We again refer to \cite{bercik2026construction}, where the primary focus is an optimized construction procedure for FSBP operators. In the present work, we adapt only minor modifications to the algorithm in \cite{huybrechs2022computation} to keep the construction as simple (and basic) as possible, since our main emphasis is on the application and performance comparison of closed and open Gaussian FSBP operators in numerical methods for hyperbolic conservation laws.} 

\begin{definition} \label{def:chebyshev}
    A set $T_k = \{u_1(t),\ldots,u_k(t)\}$ of continuous functions on $[a,b]$ is called a \emph{Chebyshev set} if
    \begin{equation} \label{eq:chebyshev_det}
        \det \begin{bmatrix}
            u_1(t_1) & \cdots & u_1(t_k)\\
            \vdots   &        & \vdots \\
            u_k(t_1) & \cdots & u_k(t_k)
        \end{bmatrix} \neq 0
    \end{equation}
    for all node distributions $a \leq t_1 < \ldots < t_k \leq b$.
\end{definition}

For simplicity, we assume that $k$ is even. 
We can always ensure this assumption by enriching the function space $(\mathcal{F}^2)'$ with an additional function, e.g., the lowest-degree polynomial not included in $(\mathcal{F}^2)'$.
We refer to \cite{huybrechs2022computation} for the case that $k$ is odd. 
The \textit{moment space} $M_k$ to a given Chebyshev set $T_k = \{u_1(t),\ldots,u_k(t)\}$ is defined as
\begin{equation} \label{eq:moment_space}
    M_k = \left\{ \mathbf{c} = [c_1,\dots,c_k]^T \in \R^k \mid c_i = \int_a^b u_i(t) \intd \sigma(t), \ i=1,\dots,k \right\}.
\end{equation}
Here, $\sigma$ can take values in all non-decreasing, right-continuous functions.
We are interested in a representation of a certain point $\mathbf{c}^*$ of this space where $\intd \sigma(t) = w(t) \intd t$. 
Since the moment space is a closed convex cone, it can be characterized as the convex conical hull of 
\begin{equation}
    \textbf{u}(t) = [u_1(t), ...,u_n(t)]
\end{equation}
for $t \in [a,b]$. Thus, we can write the point $\mathbf{c}^*$  as a linear combination of points on $\textbf{u}(t)$: 
\begin{equation}
    \mathbf{c}^* 
    		= \sum_{i=1}^{k+1} \lambda_i \textbf{u}(t_i) 
		= \left[ \sum_{i=1}^{k+1} \lambda_i u_1(t_i), \ \dots, \ \sum_{i=1}^{k+1} \lambda_i u_k(t_i) \right],
\end{equation}
where $\lambda_1, \dots, \lambda_k > 0$ and $a \leq t_1 < \ldots < t_k \leq b$.

Note that the above corresponds to a quadrature formula with $k+1$ nodes that integrates the given Chebyshev set $T_k = \{u_1(t),\ldots,u_k(t)\}$ exactly. 
If we could construct a similar quadrature formula with only $k/2$ nodes, the desired result would follow. 
In \cite{karlin1966tchebycheff}, the authors establish a theorem regarding the number of nodes required for such a quadrature formula in terms of the so-called \emph{index}.
To clarify the notion of index, consider a set of nodes $a \leq t_1 < \cdots < t_p \leq b$. 
The \textit{index} of this set is defined as
\begin{equation}\label{eq:index}
    I\Big( \{t_i\}_{i=1}^p \Big) = 
    \begin{cases} 
        p, & \text{if } t_i \in (a,b) \\
        p - \frac{1}{2}, & \text{if } t_1 = a \text{ or } t_p = b \\
        p - 1, & \text{if } t_1 = a \text{ and } t_p = b.
    \end{cases}
\end{equation}
With this, the index of a linear combination is the index of its set of nodes, and the index of a vector $\mathbf{c} \in M_k$ is the minimal index of all its possible linear combinations.

In short, \cite{karlin1966tchebycheff} then shows that every boundary point $\mathbf{c} \in \partial M_k$ can be represented as a linear combination with index $I(\mathbf{c}) < \frac{k}{2}$. 
For an interior point $\mathbf{c} \in \operatorname{int} M_k$, \cite{karlin1966tchebycheff} proves that, for a fixed $t^* \in [a,b]$, a linear combination with index $\frac{k}{2}$ exists if $t^*$ is one of the Gaussian nodes, whereas the index is $\frac{k+1}{2}$ otherwise.
Assuming that $k$ is even, a representation of $\mathbf{c}$ with index $\frac{k}{2}$ means we have either $\frac{k}{2}$ interior nodes or $\frac{k+2}{2}$ points which include the boundary points $a$ and $b$. 
These two different representations relate to Gauss-Legendre and Gauss-Lobatto quadratures.

\subsection{Construction of open GGQs}

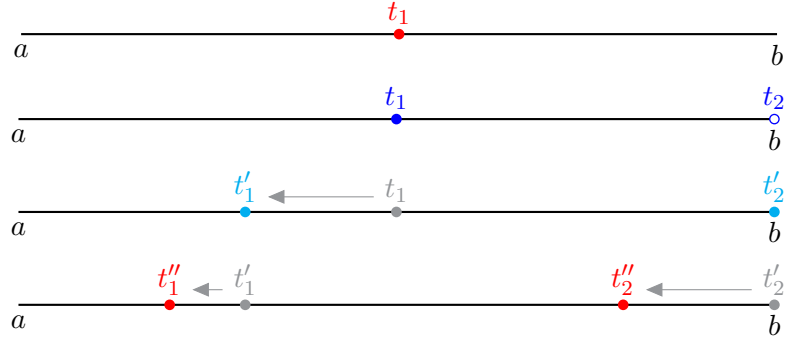
\begin{figure}[tb]
	\begin{center}
     \begin{tikzpicture}
         \coordinate (A) at (0,0);
         \coordinate (B) at (10,0);
         \draw[thick] (A) -- (B);
         \node[below] at (A) {$a$};
         \node[below] at (B) {$b$};
         \fill[red] (5,0) circle (2pt);
         \node[above] at (5,0) {\textcolor{red}{$t_1$}};
     \end{tikzpicture}
 \end{center}
 \begin{center}
     \begin{tikzpicture}
         \coordinate (A) at (0,0);
         \coordinate (B) at (10,0);
         \draw[thick] (A) -- (B);
         \node[below] at (A) {$a$};
         \node[below] at (B) {$b$};
         \fill[blue] (10,0) circle (2pt);
         \fill[blue] (5,0) circle (2pt);
         \fill[white] (10,0) circle (1.5pt);
         \node[above] at (5,0) {\textcolor{blue}{$t_1$}};
         \node[above] at (10,0) {\textcolor{blue}{$t_2$}};
     \end{tikzpicture}
 \end{center}
 \begin{center}
     \begin{tikzpicture}
         \coordinate (A) at (0,0);
         \coordinate (B) at (10,0);
         \draw[thick] (A) -- (B);
         \node[below] at (A) {$a$};
         \node[below] at (B) {$b$};
         \fill[cyan] (10,0) circle (2pt);
         \fill[gray]  (5,0) circle (2pt);
         \fill[cyan] (3,0) circle (2pt);
         \node[above] at (3,0) {\textcolor{cyan}{$t'_1$}};
         \node[above] at (5,0) {\textcolor{gray}{$t_1$}};
         \node[above] at (10,0) {\textcolor{cyan}{$t'_2$}};
         \draw[<-,gray] (3.3,0.2) -- (4.7,0.2);
     \end{tikzpicture}
 \end{center}
 \begin{center}
     \begin{tikzpicture}
         \coordinate (A) at (0,0);
         \coordinate (B) at (10,0);
         \draw[thick] (A) -- (B);
         \node[below] at (A) {$a$};
         \node[below] at (B) {$b$};
         \fill[gray] (10,0) circle (2pt);
         \fill[red]  (8,0) circle (2pt);
         \fill[gray] (3,0) circle (2pt);
         \fill[red]  (2,0) circle (2pt);
         \node[above] at (2,0) {\textcolor{red}{$t''_1$}};
         \node[above] at (3,0) {\textcolor{gray}{$t'_1$}};
         \node[above] at (8,0) {\textcolor{red}{$t''_2$}};
         \node[above] at (10,0) {\textcolor{gray}{$t'_2$}};
         \draw[<-,gray] (2.3,0.2) -- (2.7,0.2);
         \draw[<-,gray] (8.3,0.2) -- (9.7,0.2);
     \end{tikzpicture}
 	\end{center}
	\caption{Illustration of the first iteration of the algorithm in \cite{huybrechs2022computation} for computing GGQs}
	\label{fig:algo}
\end{figure}

The constructive proof in \cite{huybrechs2022computation} naturally leads to an algorithm for computing GGQs. 
We only sketch the first iteration here; the complete algorithm is given in \cite{huybrechs2022computation}. 
We illustrate the first iteration in \cref{fig:algo}.

\begin{enumerate}
    \item \textbf{Start:} Compute the 1-point quadrature for
    $u_1(t)$ and $u_2(t)$.

    \item \textbf{Step 1:} Add $b$ to the nodes and let $\xi$
    decrease from $t_1$ to $a$, and stop when $
        \sum_{i=1}^{2} \lambda_i'(\xi)\,
        u_3\!\left(t_i'(\xi)\right)
$
    integrates $u_3$ correctly.

    \item \textbf{Step 2:} Let $\xi$ decrease from $t_1''$ to $a$,
    and stop when $
        \sum_{i=1}^{2} \lambda_i''(\xi)\,
        u_4\!\left(t_i''(\xi)\right)$
    integrates $u_4$ correctly.
\end{enumerate}

Notably, we have adapted the algorithm from \cite{huybrechs2022computation} in two ways: 
First, we have increased the accuracy inside the algorithm. 
Second, we always use an orthogonal basis as the underlying input data, obtained via the Gram-Schmidt algorithm. 
Both modifications have improved the efficiency and quality of the FSBP operators.

\paragraph{Application to FSBP operators}

Since the construction of FSBP operators requires an $(\mathcal{F}^2)'$-exact quadrature, we assume that a basis of $(\mathcal{F}^2)'$ forms a Chebyshev set. This assumption allows us to apply the quadrature construction of \cite{huybrechs2022computation}.

\begin{remark}\label{rem:freeDegrees}
	Depending on the function space $\mathcal{F}$, even FSBP operators based on GGQs can fail to be interpolatory. 
	For open GGQs, this occurs whenever the dimension $M = \dim((\mathcal{F}^2)')$ exceeds twice the dimension $K = \dim(\mathcal{F})$, that is, $M > 2K$. 
	In this case, an open GGQ requires $N = \lceil M/2 \rceil$ points to be $(\mathcal{F}^2)'$-exact, and since $N \geq M/2 > K$, the resulting FSBP operator uses more points than $\dim(\mathcal{F})$. 
	These additional points introduce degrees of freedom that are not constrained by the exactness conditions and can instead be used to optimize the operator---for instance, to reduce its error on unresolved functions \cite{mattsson2014optimal,glaubitz2026summation}.
\end{remark}

\begin{remark}
	After the present work had been completed but not yet submitted, the recent work \cite{bercik2026construction} combined GGQs with an optimization approach similar to the one outlined in \cref{rem:freeDegrees}.
	A detailed comparison of the Gaussian FSBP operators developed here and the ones constructed in \cite{bercik2026construction}, as well as their application to hyperbolic conservation laws, is left for future work. 
	Here, we compare the open Gaussian FSBP operators developed here with the closed Gaussian FSBP operators constructed in \cite{hale2026summation}.
\end{remark}

%% file: 4_examples.tex
\section{Examples of generalized FSBP operators} 
\label{sec:examples}

We present a few examples of the open Gaussian FSBP operators and compare them to closed Gaussian FSBP operators and previous FSBP operators on equidistant points. 
We consider the reference interval $[x_L,x_R] = [0,1]$ and different non-polynomial function spaces.
As a sanity check, we also applied the algorithm from Section~\ref{sec:quadratures} to polynomial function spaces and recovered the Gauss--Legendre quadratures. 
For brevity, we do not report these results here.

\subsection{Exponential basis} 
\label{sub:examples_exp}

As a first example, we consider the function space $\mathcal{E}_2 = \text{span}\{1,x,e^x\}$. 
Our goal thus is to compute a quadrature for 
\begin{equation}
    (\mathcal{E}_2^2)' = \text{span} \{1,x,e^x,xe^x,e^{2x}\}.
\end{equation}

To ensure an even number of linearly independent basis functions, we enrich $(\mathcal{E}_2^2)'$ with the additional function $x^2$---the lowest-degree polynomial not already included; see the discussion in \Cref{sub:chebyshev_quad}. 

After orthogonalizing the basis functions, we obtain the following nodes and weights for the corresponding open GGQ:
\begin{table}[H]
    \centering
    \begin{tabular}{c|c|c|c}
        \textbf{x} & $0.1224$ & $0.5237$ & $0.8966$ \\
        \hline
        \textbf{w} & $0.2985$ & $0.4438$ & $0.2577$
    \end{tabular}
\end{table}
The resulting open Gaussian FSBP operator is 
\begin{equation}
    D =
    \begin{bmatrix*}[r]
       -3.4818  &  4.5475  & -1.0658 \\
       -1.3601  &  0.1420  &  1.2181 \\
        1.5399  & -5.8797  &  4.3398 \\
    \end{bmatrix*}
\end{equation}
and uses three points. 

As a comparison, using the least-squares approach from \cite{glaubitz2023summation} results in the five-point equally-spaced quadrature formula
\begin{table}[H]
    \centering
    \begin{tabular}{c|c|c|c|c|c}
        \textbf{x} & 0 & 0.25 & 0.5 & 0.75 & 1 \\
        \hline
        \textbf{w} & 0.0759763872 & 0.3620888878 & 0.1244746618 & 0.3608784643 & 0.0765815990
    \end{tabular}
\end{table} 
The resulting FSBP operator is 
\begin{equation}
    D =
    \begin{bmatrix*}[r]
        -6.580992049 & 8.594176227 & -0.461013497 & -2.536533493 & 0.984362811 \\
        -1.803298810 & 0.000000000 & 0.883330737 & 1.446533768 & -0.526565695 \\
        0.281391726 & -2.569553026 & 0.000000000 & 2.583092176 & -0.294930875 \\
        0.534020923 & -1.451385592 & -0.890963459 & 0.000000000 & 1.808328128 \\
        -0.976583554 & 2.489678844 & 0.479376527 & -8.521455370 & 6.528983553
    \end{bmatrix*}.
\end{equation}
Notably, our open Gaussian FSBP operator requires two fewer points than the above least-squares quadrature in \cite{glaubitz2023summation}.

Furthermore, in \cite[Section 3.5, Example I]{hale2026summation}, the authors obtained the following four-point closed GGQ for the same function space: 
\begin{table}[H]
    \centering
    \begin{tabular}{c|c|c|c|c}
        \textbf{x} & 0 & 0.2956452974 & 0.7423537958  & 1 \\
        \hline
        \textbf{w} & 0.0914828668 & 0.4341375639&	0.3987262252 & 0.0756533441
    \end{tabular}
\end{table} 
The resulting closed Gaussian FSBP operator is
\begin{equation}
    D =
    \begin{bmatrix*}[r]
        	-5.465504277& 7.365125959& -2.802901094&   0.903279412 \\
		-1.552003083& 0&  2.142484824& -0.590481741\\
		0.643091453& -2.332761387&            0&   1.689669934\\
		-1.092279411&  3.388486098& -8.905299859&   6.609093173
    \end{bmatrix*}.
\end{equation}
Notably, our open Gaussian FSBP operator requires one fewer point than the closed Gaussian FSBP operator in \cite{hale2026summation}.

\subsection{Hyperbolic basis}
\label{sub:examples_hyp}

We next consider the hyperbolic function space $\mathcal{H} = \text{span}\{1, x,  \sinh(x), \cosh(x)\}$. 
Such hyperbolic functions arise in dispersive wave problems, whose solutions are characterized by them. 

Our goal thus is to compute a positive quadrature formula that is exact for
\begin{equation}
    (\mathcal{H}^2)' = \text{span}\{ 1, x, \sinh(x), \cosh(x), x \sinh(x), x \cosh(x), \sinh(2x), \cosh(2x) \}.
\end{equation}
We again orthogonalize the basis functions and obtain the following quadrature formula:
\begin{table}[H]
    \centering
    \begin{tabular}{c|c|c|c|c}
        \textbf{x} & 0.0691 & 0.3296 & 0.6704 & 0.9309 \\
        \hline
        \textbf{w} & 0.1733 & 0.3267 & 0.3267 & 0.1733
    \end{tabular}
\end{table} 
The resulting open Gaussian FSBP operator is 
\begin{equation}
    D =
    \begin{bmatrix*}[r]
        -6.7200 &  9.7968 & -4.2229 &  1.1461 \\
        -1.5012 & -0.7798 &  2.9280 & -0.6471 \\
         0.6471 & -2.9280 &  0.7798 &  1.5012 \\
        -1.1461 &  4.2229 & -9.7968 &  6.7200
    \end{bmatrix*}.
\end{equation}
Using the algorithm described in \cite{hale2026summation}, we computed the following $5$-point closed GGQ quadrature:
\begin{table}[H]
    \centering
    \begin{tabular}{c|c|c|c|c|c}
        \textbf{x} & 0 & 0.1721 & 0.5 & 0.8279 & 1 \\
        \hline
        \textbf{w} & 0.0497 & 0.272 & 0.3566 & 0.272 & 0.0497
    \end{tabular}.
\end{table} 
The resulting closed Gaussian FSBP operator is 
\begin{equation}
    D =
    \begin{bmatrix*}[r]
        -10.0594 & 13.5875 & -5.3460 &   2.8097 & -0.9918 \\
         -2.4830 &       0 &  3.4892 &  -1.5196 &  0.5135 \\
          0.7451 & -2.6612 &       0 &   2.6612 & -0.7451 \\
         -0.5135 &  1.5196 & -3.4892 &        0 &  2.4830 \\
          0.9918 & -2.8097 &  5.3460 & -13.5875 & 10.0594 \\
    \end{bmatrix*}.
\end{equation}
Compared with the above closed Gaussian FSBP operator from \cite{hale2026summation}, our open Gaussian FSBP operator needs one point less.

\subsection{Bessel functions}
\label{sub:examples_Bessel}

As a final, higher-dimensional example, we consider the Bessel function space

\begin{equation}\label{eqn:example2}
    \mathcal{F} = \text{span}\{J_\nu(x)\}_{\nu=0}^{9},
\end{equation}
defined on the interval $[0,25]$, where $J_\nu(x)$ denote the Bessel functions of the first kind.
Such functions frequently arise in the modeling of radially symmetric waves and flows \cite{korenev2002bessel}.

In \cite[Section 3.5, Example II]{hale2026summation}, the same function space was considered.
Since the corresponding space $\mathcal{G} = (\mathcal{F}^2)'$ cannot easily be determined in closed form, the authors instead employ Chebfun's automatic SVD-based procedure to extract a basis of $\mathcal{G}$ numerically, reporting $\dim \mathcal{G} = 48$.
We found this value to be incorrect, however; the procedure described in \cref{rem:Bessel} below yields the correct dimension $\dim \mathcal{G} = 27$.\footnote{
We attribute the error to the severe ill-conditioning of the Vandermonde-type matrices arising from $\mathcal{G}$ on $[0,25]$, which likely misled Chebfun's automatic SVD.
}
Applying their algorithm, adapted from the GGQ method of \cite{ma1996generalized}, the authors of \cite{hale2026summation} report a positive 25-point closed GGQ that is shown in Table \ref{tab:bessel-hale}.

\begin{table}[H]
\centering
\caption{Generalized Gauss--Lobatto quadrature rule reported in \cite{hale2026summation} (25 nodes; not exact for $(\mathcal{F}^2)'$)}
\label{tab:bessel-hale}
\setlength{\tabcolsep}{10pt}%
\begin{tabular}{r S[table-format=2.10,table-column-width=2.65cm] S[table-format=1.11,table-column-width=2.65cm] | r S[table-format=2.8,table-column-width=2.25cm] S[table-format=1.11,table-column-width=2.65cm]}
	\toprule
	{$i$} & {$x_i$} & {$w_i$} & {$i$} & {$x_i$} & {$w_i$} \\
	\midrule
	1 & 0 & 0.04674109541 & 13 & 12.14501738 & 1.586926028 \\
	2 & 0.1708613339 & 0.2857413813 & 14 & 13.78368019 & 1.634993556 \\
	3 & 0.5661606569 & 0.5014293903 & 15 & 15.30246428 & 1.373778538 \\
	4 & 1.166443811 & 0.6963842628 & 16 & 16.63807436 & 1.392918441 \\
	5 & 1.954074975 & 0.8750817726 & 17 & 18.11924712 & 1.503655935 \\
	6 & 2.905578711 & 1.022716109 & 18 & 19.49873798 & 1.226087151 \\
	7 & 3.993567909 & 1.150946736 & 19 & 20.68293671 & 1.214279320 \\
	8 & 5.198004650 & 1.252002006 & 20 & 21.87902576 & 1.094192400 \\
	9 & 6.495279513 & 1.347451304 & 21 & 22.83987635 & 0.8957497727 \\
	10 & 7.897105635 & 1.446405823 & 22 & 23.70798388 & 0.7768641260 \\
	11 & 9.334639871 & 1.393772578 & 23 & 24.35339492 & 0.5561938824 \\
	12 & 10.680527 & 1.350167740 & 24 & 24.80987539 & 0.3170472565 \\
	\multicolumn{3}{c|}{\smash{\small$\vdots$}} & 25 & 25.0 & 0.05847348825 \\
	\bottomrule
\end{tabular}
\end{table}

\begin{remark}[Procedure for recovering an orthogonal basis of $(\mathcal{F}^2)'$]\label{rem:Bessel}
	Our procedure recovers the dimension and an $L^2$-orthogonal basis of a space $\mathcal{V} = \text{span}\{\phi_i\}$ given by a possibly overcomplete spanning set $\{\phi_i\}$.
	Its central object is the Gram matrix $G$ of the $L^2$ inner product,
	\begin{equation}\label{eqn:gram}
    		G_{ij} = \langle \phi_i, \phi_j \rangle_{L^2} = \int \phi_i \phi_j \intd x,
	\end{equation}
	with $\operatorname{rank} G = \dim \mathcal{V}$.
	We approximate this inner product by a high-order Gauss--Legendre rule with nodes $x_k$ and weights $w_k$, and set $V_{ki} = \phi_i(x_k)$ and $W = \diag(w_k)$.
	Because the Gauss--Legendre rule integrates the products $\phi_i \phi_j$ to high precision, the matrix $V^T W V$ reproduces $G$:
	\begin{equation}\label{eqn:gram-quad}
    		\bigl(V^T W V\bigr)_{ij} 
			= \sum_k w_k\, \phi_i(x_k)\, \phi_j(x_k) 
			\approx \int \phi_i \phi_j \intd x 
			= G_{ij}
	\end{equation}
	For readability, we treat this reproduction as exact in the following, neglecting the small quadrature error.
	The singular value decomposition $W^{1/2} V = U \Sigma Z^T$ then provides the eigendecomposition
	\begin{equation}\label{eqn:gram-svd}
    		G = V^T W V = Z\, \Sigma^2\, Z^T,
	\end{equation}
	so the squared singular values $\sigma_i^2$ are the eigenvalues of $G$, and the functions $\psi_i = \sum_j Z_{ji}\, \phi_j$ form an $L^2$-orthogonal basis of $\mathcal{V}$, with $\langle \psi_i, \psi_j \rangle_{L^2} = \sigma_i^2\, \delta_{ij}$.
	In particular, $\dim \mathcal{V}$ equals the numerical rank of $G$, the number of $\sigma_i$ above a relative threshold, here $10^{-14}\, \sigma_{\max}$.
\end{remark}

Here $\mathcal{G} = (\mathcal{F}^2)'$ is spanned by the derivatives of the $10 \cdot 11 / 2 = 55$ products $J_\mu J_\nu$ with $0 \le \mu \le \nu \le 9$, and the procedure described in \cref{rem:Bessel} returns $\dim \mathcal{G} = 27$.
We then augmented the resulting orthogonal basis with the constant function $1$ to reach an even number of $28$ linearly independent basis functions.
Notably, the algorithm of \cite{huybrechs2022computation} failed to converge for this augmented $\mathcal{G}$, so we instead used the measure-continuation method of \cite{yarvin1998generalized}.
In this way, we obtained an open GGQ with 14 points, which is reported in Table~\ref{tab:bessel-open}.

\begin{table}[tb]
  \centering
  \caption{Open quadrature rule exact for $(\mathcal{F}^2)'$ (14 nodes)}
  \label{tab:bessel-open}
  \begin{tabular}{r S S}
    \toprule
    {$i$} & {$x_i$} & {$w_i$} \\
    \midrule
    1 & 0.196619234645 & 0.499519933479 \\
    2 & 0.995125651591 & 1.07560058667 \\
    3 & 2.29881803578 & 1.50737365928 \\
    4 & 3.96568796471 & 1.80716010945 \\
    5 & 5.88062476057 & 2.00915908666 \\
    6 & 7.96129208469 & 2.14266638726 \\
    7 & 10.1496807039 & 2.22706452455 \\
    8 & 12.4027843236 & 2.27333394613 \\
    9 & 14.6851213977 & 2.28573409098 \\
    10 & 16.9622438217 & 2.26173865952 \\
    11 & 19.1928369322 & 2.18960315408 \\
    12 & 21.3142059999 & 2.03528186331 \\
    13 & 23.2082268853 & 1.70999637227 \\
    14 & 24.6010360007 & 0.97576762637 \\
    \bottomrule
  \end{tabular}
\end{table}
 
Moreover, fixing the endpoints, we obtained a closed GGQ with 15 points, which is reported in Table~\ref{tab:bessel-closed}. 

\begin{table}[tb]
  \centering
  \caption{Closed quadrature rule exact for $(\mathcal{F}^2)'$ (15 nodes)}
  \label{tab:bessel-closed}
  \begin{tabular}{r S S}
    \toprule
    {$i$} & {$x_i$} & {$w_i$} \\
    \midrule
    1 & 0 & 0.137754661776 \\
    2 & 0.492964059844 & 0.807761582166 \\
    3 & 1.56758321783 & 1.31610490922 \\
    4 & 3.07609053641 & 1.67856300834 \\
    5 & 4.88607823794 & 1.92513012997 \\
    6 & 6.89910569572 & 2.08960916967 \\
    7 & 9.04621928741 & 2.19653377265 \\
    8 & 11.2780702082 & 2.26090655343 \\
    9 & 13.5564293814 & 2.29025873684 \\
    10 & 15.8474039286 & 2.28578927178 \\
    11 & 18.1145447753 & 2.24081687754 \\
    12 & 20.3091972797 & 2.13493621685 \\
    13 & 22.3481092223 & 1.91617276092 \\
    14 & 24.054358062 & 1.4303850349 \\
    15 & 25 & 0.289277313932 \\
    \bottomrule
  \end{tabular}
\end{table}

By contrast, constructing a valid FSBP operator for the same function space using the least-squares approach of \cite{glaubitz2023construction} requires 58 equally spaced quadrature points, highlighting the efficiency of the proposed GGQs.
This 58-point rule already improves upon the 100-point equally spaced rule reported in \cite{hale2026summation}, but still highlights a substantial gap to our 14- and 15-points quadratures.

%% file: 5_numerics.tex
\section{Computational examples} 
\label{sec:numerics}

We provide computational experiments comparing closed and open Gaussian FSBP operators as well as closed and open polynomial SBP operators. 
Our computational experiments span three different one-dimensional hyperbolic conservation laws, ranging from a boundary layer problem for the inhomogeneous linear advection equation over the inviscid Burgers equation to the compressible Euler equations.
If not otherwise specified, we use a multiblock FSBP-SAT semidiscretization as in \cite{glaubitz2023summation} and an explicit strong-stability-preserving (SSP) Runge--Kutta (RK) method of third order with three stages (SSPRK(3,3)) \cite{shu1988total} for integrating in time.

\begin{remark}\label{rem:scaling}
	Polynomials are invariant under both translation and scaling, so for $\mathcal{P}_d = \text{span}\{ x^k \mid k=0,1,\dots,d \}$ a single reference operator serves every position and width, and no scaled/fixed variant is needed. 
	The exponential spaces $\mathcal{E}_d = \text{span}\{ x^k, e^x \mid k=0,1,\dots,d-1 \}$ are also translation invariant, as $e^{x+c}=e^{c}\,e^{x}$. 
	They are not scale invariant, however, since $e^{sx}\notin\mathcal{E}_d$ for $s\neq1$.
	Hence, their scale must be matched to the blocks through the block size.
\end{remark}

\subsection{A boundary layer problem for the inhomogeneous linear advection} 
\label{sub:linear}

Consider the  following inhomogeneous linear advection problem from \cite[Section 7.2]{glaubitz2023summation}: 
\begin{equation}\label{eq:linear_inhom} 
\begin{aligned} 
	\partial_t u + \partial_x u & = 2u, \quad && 0 < x < \pi, \\ 
	u(x,0) & = 1, \quad && 0 \leq x \leq \pi, \\ 
	u(0,t) & = 1, \quad && t \geq 0,
\end{aligned} 
\end{equation} 
with exact steady state solution $u(x) = e^{2x}$. 
The steady state solution is likely to be better approximated using an exponential than a polynomial approximation space.  
Thus, for this example, we compared different exponential operators:  
First, we used the least-squares FSBP operator from \cite{glaubitz2023summation}, which uses five equidistant nodes and is exact for the space $\mathcal{E}_2 = \{1, x, e^x\}$. Next, we chose the closed Gaussian FSBP operator from \cite{hale2026summation}, which is exact for the same space but uses only four nodes. 
And finally, we compare these two FSBP operators with our open Gaussian FSBP operator, which also uses four nodes but is exact for the larger space $\mathcal{E}_3 = \{1, x, x^2, e^x\}$.

\begin{figure}[tb]
    \begin{subfigure}{.49\textwidth}
        \centering
        \includegraphics[width=\linewidth]{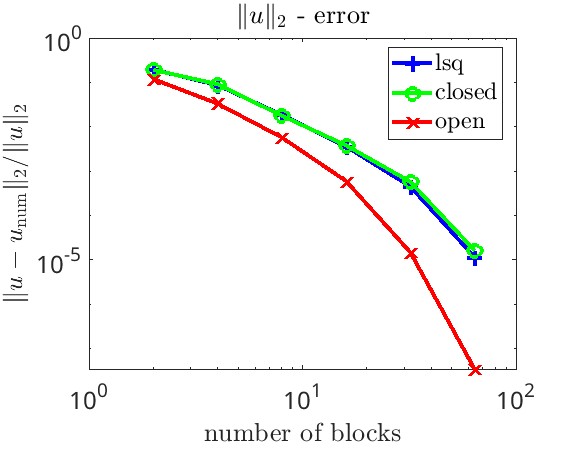}
        \caption{$2$-norm errors depending on the number of blocks}
    \end{subfigure}
    \begin{subfigure}{.49\textwidth}
        \centering
        \includegraphics[width=\linewidth]{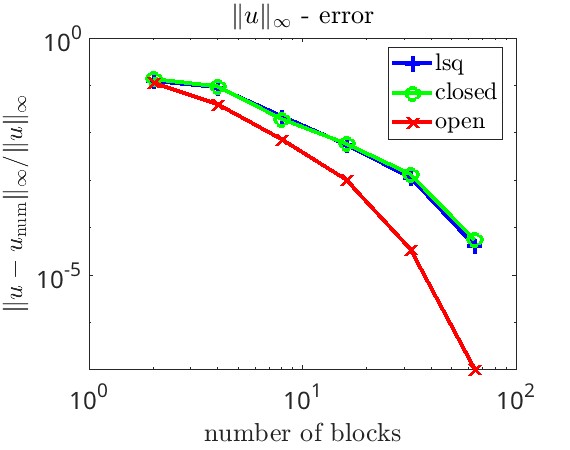}
        \caption{$\infty$-norm errors depending on the number of blocks}
    \end{subfigure}
    \caption{
		Errors of the numerical solutions for the inhomogeneous linear advection equation \cref{eq:linear_inhom}.
		We compare the least-squares FSBP operator from \cite{glaubitz2023summation} (denoted ``lsq''), the closed Gaussian FSBP operator from \cite{hale2026summation} (denoted ``closed''), and our open Gaussian FSBP operator (denoted ``open'').
	}
    \label{fig:inhomog_adv}
\end{figure}

\cref{fig:inhomog_adv} displays the relative $\Vert \cdot \Vert_2$- and $\Vert \cdot \Vert _{\infty}$-errors of the numerical solutions obtained using a multi-block FSBP-SAT method with the different FSBP operators at time $T = 3.5$ for different numbers of blocks.
Using $2^k$ equisized blocks for $k = 1,...,6$ and computing a new operator on $[0, \frac{\pi}{2^k}]$ for each refinement level, we observe that the $4 \times 4$ and $5 \times 5$ FSBP operators with the least-squares quadrature and closed GGQ lead to rather similar results. 
Using our $4 \times 4$ open Gaussian FSBP operator instead, the errors get significantly smaller. 
This can be explained by the fact that this operator is exact for a larger approximation space using the same number of points as the closed Gaussian FSBP operator.

\subsection{Inviscid Burgers' equation}
\label{sub:burgers}

Following \cite[Section 7.3]{glaubitz2023summation}, we consider the inviscid Burgers' equation given by
\begin{equation}\label{eq:burgers}
    \begin{aligned}
        \partial_tu + \partial_x \Big(\frac{u^2}{2}\Big) &= 0, \quad && 0 < x < 1 \\
        u(x,0) &= 1+ \frac{1}{2}\sin(4\pi x)^3 + \frac{1}{4} \cos(4\pi x )^5, \quad && 0 \leq x \leq 1 \\
        u(0,t) &= u(1,t), \quad && t \geq 0.
    \end{aligned}
\end{equation}

We compare the same four-point closed and open FSBP operators, respectively exact for $\mathcal{E}_2 = \{1, x, e^x\}$ and $\mathcal{E}_3 = \{1, x, x^2, e^x\}$, as in \Cref{sub:linear}. 
Furthermore, we combine the FSBP-SAT semidiscretization with a skew-symmetric formulation of \cref{eq:burgers}; see \cite[Section 7.3]{glaubitz2023summation}.

\begin{figure} 
    \begin{subfigure}{.49\textwidth}
        \centering
        \includegraphics[width=\linewidth]{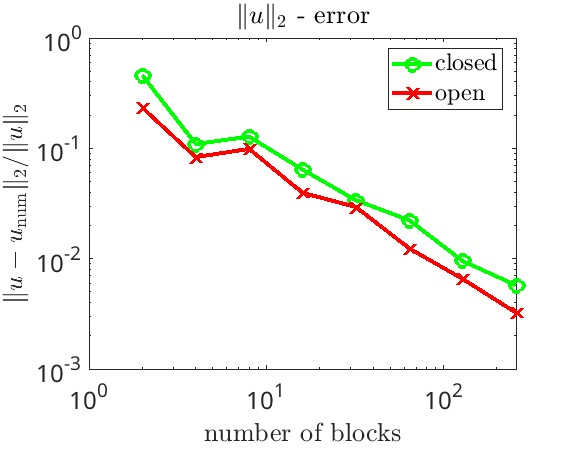}
        \caption{$2$-norm errors, reference block $[0,1]$}
    \end{subfigure}
    \begin{subfigure}{.49\textwidth}
        \centering
        \includegraphics[width=\linewidth]{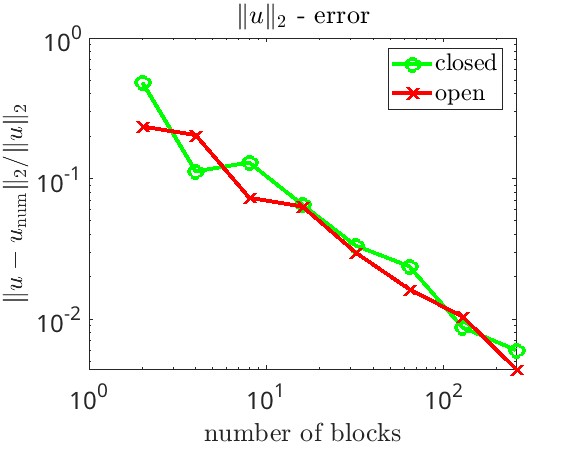}
        \caption{$2$-norm errors, reference blocks $[0, \frac{1}{2^k}]$}
    \end{subfigure}
    \caption{
		Errors of the numerical solutions for the inviscid Burgers equation \cref{eq:burgers}.
		We compare the closed Gaussian FSBP operator from \cite{hale2026summation} (denoted ``closed'') with our open Gaussian FSBP operator (denoted ``open'') for two different reference blocks.
	}
    \label{fig:burgers}
\end{figure}

\cref{fig:burgers} displays the relative $\Vert \cdot \Vert_2$-errors of the numerical solutions obtained with the closed Gaussian FSBP operator from \cite{hale2026summation} and our open Gaussian FSBP operator at time $T = 0.1$ for different numbers of blocks. 
We also compared two different kinds of reference blocks for the $2^k$ blocks, where $k=1,\dots,8$. 
On the one hand, we used $[0,1]$ as a reference block for all block sizes and scaled and shifted the operator and its nodes by the block size. 
On the other hand, we used $[0, \frac{1}{2^k}]$ as the reference block for each refinement level, so we do not need to scale the operator. 
Also see the discussion in \cref{rem:scaling}.
Comparing the two versions, we notice that the convergence rates appear to be similar. 
However, the individual errors for the open quadrature are slightly smaller for the $[0,1]$ reference block.

\subsection{One-dimensional compressible Euler equations}
\label{sub:Euler}

We consider the one-dimensional compressible Euler equations
\begin{equation}\label{eq:Euler}
    \frac{\partial}{\partial t}
    \begin{pmatrix}
        \rho \\
        \rho v \\
        \rho E
    \end{pmatrix}
    +
    \frac{\partial}{\partial x}
    \begin{pmatrix}
        \rho v \\
        \rho v^2 + p \\
        (\rho E + p)v
    \end{pmatrix}
    = \mathbf{0},
\end{equation}
where $\rho$ denotes the density, $v$ the velocity, $E$ the specific total energy, and the pressure is given by the ideal gas law
\begin{equation}
	p = (\gamma-1)\left(\rho E - \frac{1}{2}\rho v^2\right).
\end{equation}
Throughout the experiments, we choose the ratio of specific heats $\gamma = 1.4$.

The computational domain is $\Omega=[-1,1]$ with periodic boundary conditions, and the solution is evolved until the final time $T=2$.
As an exact smooth solution of the homogeneous Euler equations, we consider the advecting density wave
\begin{equation}\label{eq:initial_condition_Euler}
\begin{aligned}
    \rho(x,t) &= 1 + \frac{1}{2}\sin\!\bigl(2\pi(x-0.1t)\bigr),\\
    (\rho v)(x,t) &= 0.1\,\rho(x,t),\\
    (\rho E)(x,t) &= 50 + 0.005\,\rho(x,t).
\end{aligned}
\end{equation}
This solution corresponds to a constant velocity $v=0.1$ and constant pressure $p=20$, and is therefore an exact source-free solution of the full nonlinear Euler equations. 
The initial condition is obtained by evaluating \cref{eq:initial_condition_Euler} at $t=0$.

All computations are carried out with the Trixi.jl framework
\cite{ranocha2022adaptive}. 
For all discretizations, we employ the HLLC numerical flux \cite{toro1994restoration} at the element interfaces and discretize the volume terms in weak form using the SBP derivative operators. 
Time integration uses the explicit five-stage, fourth-order low-storage Carpenter--Kennedy RK method \cite{carpenterFourthorder2NstorageRungeKutta1994} (the \texttt{CarpenterKennedy2N54} integrator of OrdinaryDiffEq.jl \cite{rackauckas2017differentialequations}) with an adaptive time step.
We use uniform meshes of $2^\ell$ elements at refinement levels $\ell=2,\dots,8$, with element width $h=2/2^\ell$. 

We compare diagonal-norm SBP discretizations built from the polynomial spaces $\mathcal{P}_d$ of degree at most $d$ and the exponential spaces $\mathcal{E}_2=\operatorname{span}\{1,x,e^x\}$ and $\mathcal{E}_3=\operatorname{span}\{1,x,x^2,e^x\}$ on an interval $[0,b]$.
The closed Gaussian FSBP operators include the two-element boundaries in their nodal set, so interface values are read off directly. 
Our open Gaussian FSBP operators use interior nodes only, with interface values obtained from an extrapolation
operator that is exact on the same function space; see \Cref{sub:const_B_FSBP}. 
For each exponential space we use two operators: a \emph{scaled} one with $b=h$, so that on every
element the space is $\{1,x,e^x\}$, independent of the mesh, and a \emph{fixed} one with $b=1$ regardless of $h$. 
The compared operators are listed in Table~\ref{tab:euler_operators}.

\begin{table}[!ht]
\centering
\caption{Operators compared in the convergence study, by function space (rows) and node distribution (columns). 
Each exponential entry denotes two operators, scaled ($b=h$) and fixed ($b=1$).}
\label{tab:euler_operators}
\begin{tabular}{lll}
\toprule
 & closed & open \\
\midrule
$\mathcal{P}_2$ & 3 Gauss--Lobatto nodes & 3 Gauss nodes \\
$\mathcal{P}_3$ & 4 Gauss--Lobatto nodes & 4 Gauss nodes \\
$\mathcal{E}_2$ & 4 closed GGQ nodes (scaled and fixed) & 3 open GGQ nodes (scaled and fixed) \\
$\mathcal{E}_3$ & --- & 4 nodes (scaled and fixed) \\
\bottomrule
\end{tabular}
\end{table}

We report the discrete $L^2$- and $L^\infty$-errors for $\rho$ at the final time.
The behavior for the other conserved variables, namely $\rho v$ and $\rho E$, is qualitatively the same, but omitted for the sake of space.

The discrete $L^2$-error is computed on each element's own quadrature nodes,
\begin{equation}
\|u-u_h\|_{L^2}
=
\left(
\frac{1}{2}
\sum_{K}
\sum_{i}
\omega_i
\bigl(u(x_i)-u_h(x_i)\bigr)^2
\right)^{1/2},
\end{equation}
where the $\omega_i$ are the quadrature weights, summing to the domain length $2$ (hence the prefactor $1/2$). 
The $L^\infty$-error is the maximum absolute nodal error over all elements. 
Experimental orders of convergence (EOCs) are computed as
\begin{equation}
\mathrm{EOC}
=
\log_2
\left(
\frac{\mathrm{error}_{\mathrm{coarse}}}
     {\mathrm{error}_{\mathrm{fine}}}
\right)
\end{equation}
between successive refinement levels, which differ by a factor of two in the number of cells (so the base-two logarithm directly yields the order).

\subsubsection*{Results}

\begin{figure}[p]
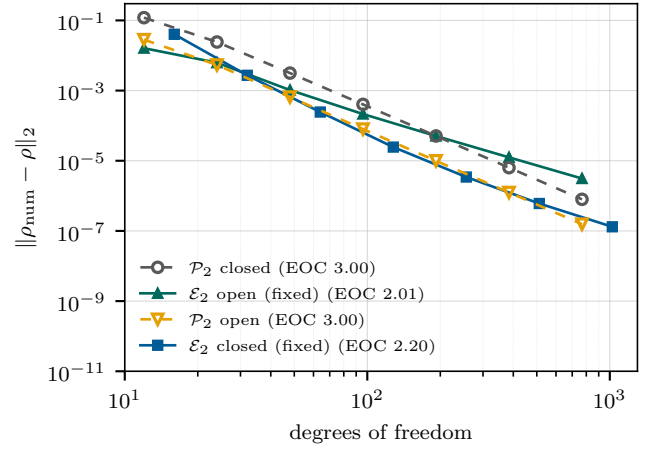
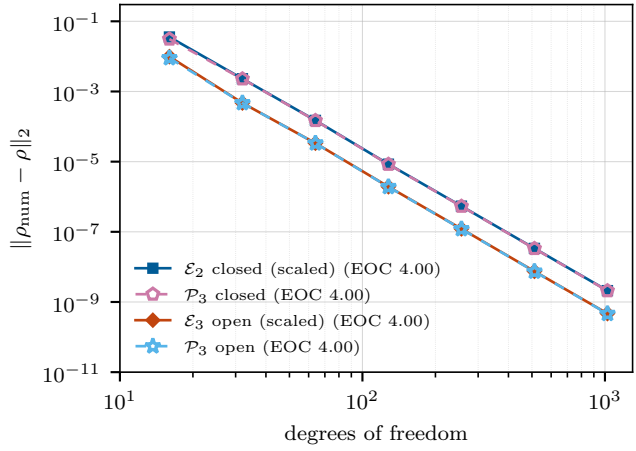
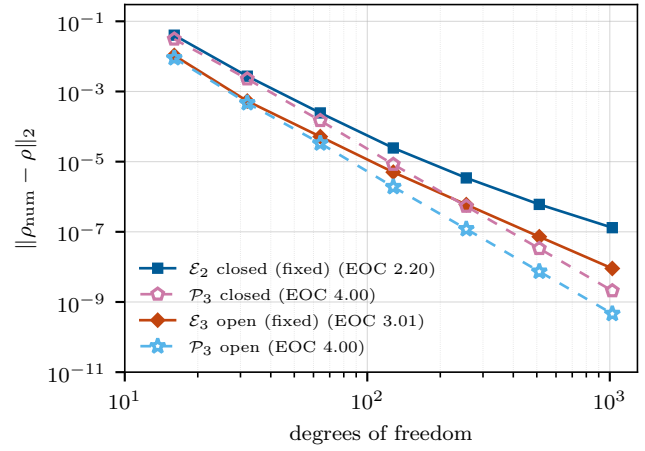
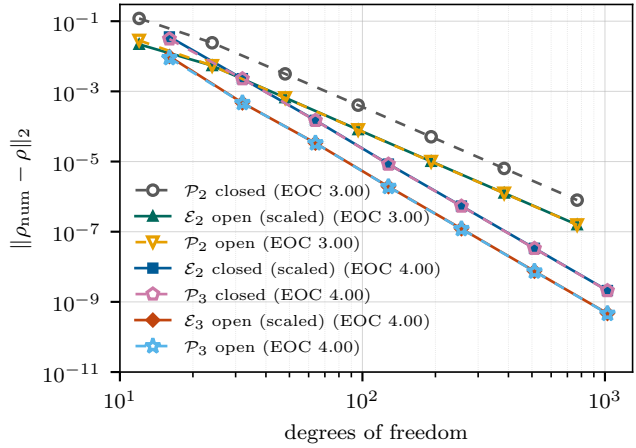
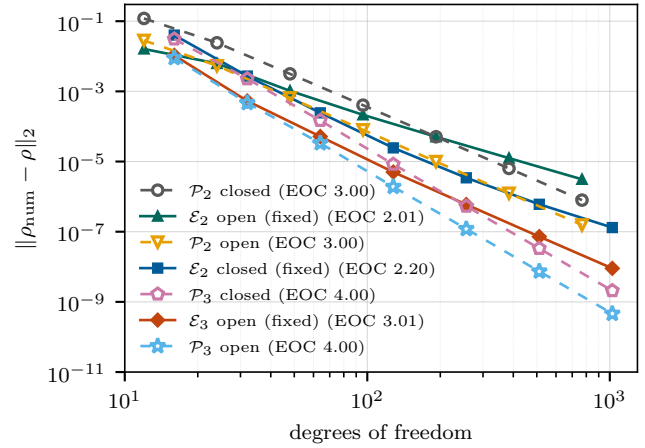

    \centering
    \begin{subfigure}{.49\textwidth}
        \centering
        \input{error_vs_nodes.pgf}
        \caption{$\mathcal{P}_2,\mathcal{E}_2$ (scaled).}
        \label{fig:conv-l2-p2e2-scaled}
    \end{subfigure}
    \hfill
    \begin{subfigure}{.49\textwidth}
        \centering
        \input{error_vs_nodes_fixed.pgf}
        \caption{$\mathcal{P}_2,\mathcal{E}_2$ (fixed).}
        \label{fig:conv-l2-p2e2-fixed}
    \end{subfigure}

    \vspace{0.2em}

    \begin{subfigure}{.49\textwidth}
        \centering
        \input{error_vs_nodes2.pgf}
        \caption{Four-node operators (scaled).}
        \label{fig:conv-l2-4node-scaled}
    \end{subfigure}
    \hfill
    \begin{subfigure}{.49\textwidth}
        \centering
        \input{error_vs_nodes2_fixed.pgf}
        \caption{Four-node operators (fixed).}
        \label{fig:conv-l2-4node-fixed}
    \end{subfigure}

    \vspace{0.2em}

    \begin{subfigure}{.49\textwidth}
        \centering
        \input{error_vs_nodes3.pgf}
        \caption{All operators (scaled).}
        \label{fig:conv-l2-all-scaled}
    \end{subfigure}
    \hfill
    \begin{subfigure}{.49\textwidth}
        \centering
        \input{error_vs_nodes3_fixed.pgf}
        \caption{All operators (fixed).}
        \label{fig:conv-l2-all-fixed}
    \end{subfigure}

    \caption{$L^2$-error of $\rho$ at $T=2$ versus degrees of freedom;
    legend EOCs from the two finest levels.}
    \label{fig:conv-l2}
\end{figure}

The $L^2$- and $L^\infty$-errors are shown in Figures~\ref{fig:conv-l2} and \ref{fig:conv-linf}, whose legends report each operator's EOC from the two finest levels. 
In each figure, the scaled operators are in the left column and the fixed operators in the right. 
The three rows compare operators of the same function-space dimension ($\mathcal{P}_2$ and $\mathcal{E}_2$), operators with the same number of nodes (four), and all operators together. 
We report the density $\rho$; the other conserved variables behave similarly.

For the scaled operators (left columns), every operator converges with EOC equal to its number of nodes. 
The scaled exponential operators attain the EOC of the polynomial operator with the same number of nodes: 
The open $\mathcal{E}_2$ operator (three nodes) reaches EOC~$3$, like $\mathcal{P}_2$, and the closed $\mathcal{E}_2$ and open $\mathcal{E}_3$ operators (four nodes) reach EOC~$4$, like $\mathcal{P}_3$. 

Moreover, the open $\mathcal{E}_2$ curve nearly coincides with the open $\mathcal{P}_2$ curve (Figures~\ref{fig:conv-l2-p2e2-scaled},~\ref{fig:conv-linf-p2e2-scaled}). 
In the four-node subplots (Figures~\ref{fig:conv-l2-4node-scaled},~\ref{fig:conv-linf-4node-scaled}), the closed $\mathcal{E}_2$ curve nearly coincides with the closed $\mathcal{P}_3$ curve and the open $\mathcal{E}_3$ curve with the open $\mathcal{P}_3$ curve. 
The open Gaussian FSBP operators are slightly more accurate than the closed Gaussian FSBP operators at equal degrees of freedom.

Fixing the exponential scale to $b=1$ (right columns) reduces the EOC: the error
curves bend toward a shallower slope, and the open and closed $\mathcal{E}_2$
operators reach EOC~$\approx 2$ while the open $\mathcal{E}_3$ operator reaches
EOC~$\approx 3$. As these slopes are shallower than those of the scaled and
polynomial operators, the fixed exponentials lose accuracy under refinement: the
fixed open $\mathcal{E}_2$ operator is the most accurate at the lowest degrees of
freedom (Figure~\ref{fig:conv-l2-all-fixed}) but the least accurate at the
highest.

\providecommand{\mathdefault}[1]{#1}    
\newcommand{\pgfdir}{figures/}

\begin{figure}[p]
    \centering
    \begin{subfigure}{.49\textwidth}
        \centering
        \input{error_vs_nodes_linf.pgf}
        \caption{$\mathcal{P}_2,\mathcal{E}_2$ (scaled).}
        \label{fig:conv-linf-p2e2-scaled}
    \end{subfigure}
    \hfill
    \begin{subfigure}{.49\textwidth}
        \centering
        \input{error_vs_nodes_fixed_linf.pgf}
        \caption{$\mathcal{P}_2,\mathcal{E}_2$ (fixed).}
        \label{fig:conv-linf-p2e2-fixed}
    \end{subfigure}

    \vspace{0.2em}

    \begin{subfigure}{.49\textwidth}
        \centering
        \input{error_vs_nodes2_linf.pgf}
        \caption{Four-node operators (scaled).}
        \label{fig:conv-linf-4node-scaled}
    \end{subfigure}
    \hfill
    \begin{subfigure}{.49\textwidth}
        \centering
        \input{error_vs_nodes2_fixed_linf.pgf}
        \caption{Four-node operators (fixed).}
        \label{fig:conv-linf-4node-fixed}
    \end{subfigure}

    \vspace{0.2em}

    \begin{subfigure}{.49\textwidth}
        \centering
        \input{error_vs_nodes3_linf.pgf}
        \caption{All operators (scaled).}
        \label{fig:conv-linf-all-scaled}
    \end{subfigure}
    \hfill
    \begin{subfigure}{.49\textwidth}
        \centering
        \input{error_vs_nodes3_fixed_linf.pgf}
        \caption{All operators (fixed).}
        \label{fig:conv-linf-all-fixed}
    \end{subfigure}

    \caption{$L^\infty$-error of $\rho$ at $T=2$ versus degrees of freedom;
    legend EOCs from the two finest levels.}
    \label{fig:conv-linf}
\end{figure}

%% file: 6_summary.tex
\section{Summary}
\label{sec:summary}

We have investigated the use of closed and open GGQs to construct efficient FSBP operators and operationalized them for hyperbolic conservation laws.
Since open GGQs exclude the interval endpoints, we first extended the FSBP framework in \Cref{sec:FSBP} to grids that need not contain one or both boundary nodes, with a boundary operator $B = \mathbf{t}_R \mathbf{t}_R^T - \mathbf{t}_L \mathbf{t}_L^T$ assembled from $\mathcal{F}$-exact extrapolation vectors; see \cref{alg:FSBP_construction}.

The efficiency gain of open over closed GGQs is, at its core, a counting statement: an open GGQ uses $m$ interior nodes to be exact for a $2m$-dimensional function system, whereas a closed GGQ needs $m+1$ nodes, two of which are pinned to the endpoints.
Applied to $(\mathcal{F}^2)'$ with $M = \dim (\mathcal{F}^2)'$, an open GGQ thus requires $N = \ceil{M/2}$ points and a closed GGQ exactly one point more.
Closed GGQs therefore remain optimal within the class of closed FSBP operators, while open GGQs extend this optimality to the open setting whenever an open formulation is admissible.

In \Cref{sec:numerics}, we applied the resulting operators to a boundary-layer problem for the inhomogeneous linear advection equation, Burgers' equation, and the one-dimensional compressible Euler equations.
At an equal number of nodes, the open construction buys exactness in a strictly larger function space, which translates into markedly smaller errors for the boundary-layer problem, whereas open and closed operators converged at comparable rates for Burgers' equation.
For the Euler equations, every operator except the fixed-scale exponentials converged at an order equal to its number of nodes.

The gains for open Gaussian FSBP operators come at a price, however.
Boundary and interface data must be extrapolated rather than read off directly.
The construction moreover presupposes that $(\mathcal{F}^2)'$ forms a Chebyshev set; and even when it holds, the algorithm of \cite{huybrechs2022computation} failed to converge for the Bessel space, where we resorted to the measure-continuation method of \cite{yarvin1998generalized}.

Future work will combine GGQ-based FSBP operators with additional optimization strategies, exploiting the surplus degrees of freedom that arise in the non-interpolatory case (\cref{rem:freeDegrees}) to further reduce the operator error on unresolved functions \cite{mattsson2014optimal,glaubitz2026summation}.
A detailed comparison with \cite{bercik2026construction} and the extension of open GGQ-based operators to the multi-dimensional FSBP setting \cite{glaubitz2023multi} are planned as well.

%% file: acknowledgements.tex
JG acknowledges support from the Swedish Research Council (VR) Starting Grant \#2025-05370, the Zenith Career Development Grant \#26.07, and the National Academic Infrastructure for Supercomputing in Sweden (NAISS) grants \#2025/22-1599 and \#2024/22-1207.
P\"O was supported by the German Research Foundation (DFG) within the framework of the Priority Program SPP 2410 under project no. 525866748 (OE 661/5-1) and by DFG project no. 520756621 (OE 661/4-1).

\subsection*{Declaration} The authors disclose that generative AI assistants (ChatGPT and Claude Code) were used during the preparation of this work to assist
with code development, data analysis, and text formatting. All AI-assisted code and
analyses were tested, reviewed, and verified by the authors, and all text was reviewed
and edited by the authors.